\documentclass[12pt]{amsart}
\usepackage{amsfonts, amssymb,amsmath, amsthm, amsxtra, latexsym, mathrsfs, amscd}
\usepackage{amsmath, amssymb, amsthm, times, float}
\usepackage{graphics}
\usepackage[latin1]{inputenc} %

\newtheorem{theo+}{Theorem}[section]
\newtheorem{prop+}[theo+]{Proposition}
\newtheorem{coro+}[theo+]{Corollary}
\newtheorem{lemm+} [theo+]{Lemma}
\newtheorem{deep+}  [theo+]  {Deep Result}
\newtheorem{fact+}  [theo+]  {Fact}
\theoremstyle{definition}
\newtheorem{exam+}  [theo+]  {Example}
\newtheorem{rema+}  [theo+]  {Remark}
\newtheorem{defi+}  [theo+]  {Definition}
\newtheorem{xca+}[theo+]{Exercise}

\numberwithin{equation}{section}

\usepackage{color}

\def\draft{\centerline{(Draft {\the \day}/{\the\month} \the \year.)}}
\bigskip

\def\refn#1.#2{\expandafter\def\csname#1\endcsname{[#2]}}
\def\refnr#1.{\csname#1\endcsname}

\def\fa{\mathfrak a}

\def\fg{\mathfrak g}
\def\fk{\mathfrak k}
\def\fh{\mathfrak h}
\def\fl{\mathfrak l}

\def\fp{\mathfrak p}
\def\fq{\mathfrak q}

\def\fsl{\mathfrak sl}
\def\fv{\mathfrak v}
\def\fu{\mathfrak u}
\def\fsu{\mathfrak su}
\def\a{\alpha}

\def\Claminv2{|C(\Lambda)|^{-2}}

\def\Ga{\Gamma}

\def\de{d\varepsilon}

\def\U#1#2{U_{#1}^{(#2)}}

\def\Aa2D{A^{\a,2}(D)}
\def\bAa2D{\overline{A^{\a,2}(D)}}
\def\Ab2D{A^{\beta,2}(D)}
\def\bAb2D{\overline{A^{\beta,2}(D)}}

\def\Norm#1_#2{\Vert#1\Vert_{#2}}

\def\phipl12{\phi_{p_{l_1}, p_{l_2}}}
\def\phip01{\phi_{p_{0}, p_{0}}}

\def\a{\alpha}

\def\Claminv2{|C(\Lambda)|^{-2}}

\def\Ga{\Gamma}

\def\ad{\operatorname{ad}}
\def\Ad{\operatorname{Ad}}

\def\det{\operatorname{det}}
\def\diag{\operatorname{diag}}

\def\diag{\operatorname{diag}}

\def\tr{\operatorname{tr}}

\def\End{\operatorname{End}}

\def\de{d\varepsilon}

\def\U#1#2{U_{#1}^{(#2)}}

\def\Aa2D{A^{\a,2}(D)}
\def\bAa2D{\overline{A^{\a,2}(D)}}
\def\Ab2D{A^{\beta,2}(D)}
\def\bAb2D{\overline{A^{\beta,2}(D)}}

\def\phipl12{\phi_{p_{l_1}, p_{l_2}}}
\def\phip01{\phi_{p_{0}, p_{0}}}

\def\bc{\mathbb C}
\def\barc{\overline{\mathbb C}}
\def\br{\mathbb R}
\def\bn{\mathbb N}
\def\bz{\mathbb Z}
\def\bq{\mathbb Q}
\def\bh{\mathbb H}

\def\alg/{algebra}

\def\Alg/{Algebra}

\def\alt/{alternative} 
\def\anal/{analytic}
\def\analfunc/{\anal/\ \func/}
\def\Ans/{\it Answer. \normal}
\def\ass/{associative}
\def\nass/{non-\ass/}
\def\autom/{automorphism}
\def\homom/{homomorphism}
\def\isom/{isomorphism}
\def\bdd/{bounded}
\def\Bdd/{Bounded}
\def\bddsymdom/{bounded \sym/ \dom/}
\def\Cartdom/{Cartan \dom/}
\def\bdry/{boundary}
\def\bsd/{\bdd/ \symdom/}
\def\bv/{boundary value}
\def\cf/{{\it cf}\.}
\def\Cf/{{\it Cf}\.}
\def\charr/{character}
\def\coeff/{coefficient}
\def\comm/{commutative}
\def\cpct/{compact}
\def\compl/{complex}
\def\comp/{complex}
\def\Comp/{Complex}
\def\conf/{conformal}
\def\conj/{conjugate}
\def\conn/{connect}
\def\cont/{continuous}
\def\conv/{converge} 
\def\convc/{convergence}
\def\convt/{convergent}
\def\convx/{convex}
\def\coord/{coordinate}
\def\lcoord/{local coordinate}
\def\Corr/{Corresponding}
\def\corr/{corresponding}
\def\corrd/{correspond}
\def\cov/{covariant}
\def\decomp/{decomposition}
\def\deco/{decompose}
\def\diff/{different} 
\def\Diff/{Different} 
\def\dimn/{dimension} 
\def\distr/{distribution} 
\def\div/{diverge} 
\def\dom/{domain}
\def\eg/{\hbox{\it e.g}\.}
\def\eigenf/{eigen\-\func/}
\def\eigensp/{eigen\-space}
\def\eigenv/{eigen\-value}
\def\eq/{equation}
\def\equa/{equation}
\def\de/{\diff/ial \equa/}
\def\do/{\diff/ial operator}
\def\ode/{ordinary \de/}
\def\pde/{partial \de/}
\def\pdo/{partial \diff/ial operator}
\def\psdo/{pseudo \diff/ial operator}
\def\fin/{finite}
\def\Ex/{\it Example.\ \normal}
\def\Exnr#1/{\it Example #1.\ \normal}
\def\foll/{follow}
\def\follg/{following}
\def\Follg/{Following}
\def\func/{function}
\def\Func/{Function}
\def\Fonc/{Fonc\-tion}
\def\fonc/{fonc\-tion}
\def\Funk/{Funk\-tion}
\def\funk/{Funk\-tion}
\def\gen/{general}
\def\har/{harmonic}
\def\Hint/{\it Hint. \normal}
\def\hist/{historic}
\def\histcl/{historical}
\def\hol/{holo\-morphic}
\def\homog/{ho\-mo\-ge\-ne\-ous}
\def\hyp/{hyper\-bolic}
\def\hyperg/{hyper\-geometric}
\def\ie/{\hbox{\it i.e.}}
\def\iff/{if and only if}
\def\ineq/{inequality}
\def\infra/{{\it inf\-ra}}
\def\ultra/{{\it ult\-ra}}
\def\Inpart/{In particular}
\def\inpart/{in particular}
\def\instof/{instead of}
\def\interps/{interpolation space}
\def\interp/{interpolation}
\def\Interp/{Interpolation}
\def\interpr/{Interpretation}
\def\Intr/{Introduction}
\def\intv/{interval}
\def\inv/{invariant}
\def\invc/{invariance}
\def\Iowords/{In other words}
\def\iowords/{in other words}
\def\ipr/{inner product}
\def\irred/{irreducible}
\def\lb/{line bundle}
\def\lin/{linear}
\def\lhs/{left hand side}
\def\rhs/{right hand side}
\def\loc/{local}
\def\math/{mathematic}
\def\mathcn/{\math/ian}
\def\manif/{manifold}
\def\meas/{measure}
\def\measl/{measurable}
\def\mero/{mero\-morphic}
\def\mon/{monomial}
\def\monog/{monogenic}
\def\mult/{multiple}
\def\multy/{multiply}
\def\multn/{multiplication}
\def\nas/{necessary and sufficient}
\def\nbd/{neighborhood}
\def\neg/{negative}
\def\nondeg/{nondegenerate}
\def\Oohand/{On the other hand}
\def\oohand/{on the other hand}
\def\Oonhand/{On the one hand}
\def\oonhand/{on the one hand}
\def\oper/{operator}
\def\orth/{ortho\-gonal}
\def\orthon/{ortho\-normal}
\def\otoh/{on the other hand}
\def\quat/{quaternion}
\def\pp/{\hbox{a. e.}}
\def\psh/{plurisubharmonic}
\def\pol/{polynomial}
\def\pot/{potential}
\def\pos/{positive}
\def\princ/{principle}
\def\prob/{probability}
\def\proj/{projective}
\def\projn/{projection}
\def\Proof/{\it Proof:\normal}
\def\Rem/{\it Remark\normal}
\def\Remnr#1/{\it Remark\ \normal #1. }
\def\rep/{representation}
\def\reps/{representations}
\def\meta/{metaplectic representation}
\def\repr/{reproducing}
\def\reprker/{reproducing kernel}
\def\resp/{respective} 
\def\resply/{respectively}
\def\restr/{restriction}
\def\sa/{self-adjoint}
\def\st/{such that}

\def\sol/{solution}
\def\ru/{space}
\def\sph/{spherical}
\def\ssp/{sub\ru/}
\def\sym/{symmetric}
\def\Sym/{Symmetric}
\def\symb/{symbol}
\def\symbc/{symbolic}
\def\symdom/{\sym/ domain}
\def\symp/{symplectic}
\def\Theor#1/{\fet Theorem #1.\ \normal}
\def\Lem#1/{\fet Lemma #1.\ \normal}
\def\Lemma/{\fet Lemma.\ \normal}
\def\topl/{topology}
\def\topll/{topological}
\def\transf/{transform}
\def\transl/{translation}
\def\transfn/{transformation}
\def\transv/{transvectant}
\def\trig/{trigonometric}
\def\tril/{trilinear}
\def\trilf/{trilinear form}
\def\uhp/{upper halfplane}
\def\uhs/{upper halfspace}
\def\vb/{vector bundle}
\def\vf/{vector field}
\def\vsp/{vector space}
\def\wrt/{with respect to}
\def\Wlog/{Without loss of generality}
\def\a{\alpha}

\def\E{\Varepsilon}

\def\Ab/{Abel}
\def\Ban/{Banach}
\def\Bansp/{\Ban/ space}
\def\Belt/{Bel\-tra\-mi}
\def\Berg/{Berg\-man}
\def\Bern/{Ber\-nou\-lli}
\def\Berz/{Berezin}
\def\Bess/{Bessel}
\def\Cart/{Car\-tan}
\def\Cay/{Cay\-ley}
\def\CG/{Clebsch-Gordan}
\def\Cl/{Clifford}
\def\CR/{Cauchy-Rie\-mann}
\def\Dir/{Dirichlet}
\def\Eucl/{Euclide}
\def\Eucln/{Euclidean}
\def\F/{Fourier}
\def\Hank/{Hankel}
\def\Hankf/{\Hank/ form}
\def\Herm/{Hermite}
\def\Hilb/{Hilbert}
\def\Hilbs/{Hilbert space}
\def\Hilbsp/{Hilbert space}
\def\HS/{Hilbert-Schmidt}
\def\Lag/{La\-grange}
\def\Lap/{La\-place}
\def\LapBelt/{\Lap/-\Belt/}
\def\Leb/{Lebesgue}
\def\Marc/{Mar\-cin\-kie\-wicz}
\def\Moeb/{Moebius}
\def\Moebt/{Moebius transformation}
\def\Moebtransfn/{Moebius transformation}
\def\Pla/{Plan\-che\-rel}
\def\Poin/{Poin\-car\'e}
\def\Riem/{Rie\-mann}
\def\Riemn/{\Riem/ian}
\def\psRiemn/{pseudo-\Riem/ian}
\def\Riems/{Rie\-mann surface}
\def\Schroe/{Schr\"odinger}
\def\Weier/{Weier\-strass}
\def\anal/{analytic}
\def\bsd/{bounded symmetric domain  }
\def\bdd/{bounded}
\def\calc/{calculation}\def\conj{conjugate}
\def\calci/{calculating}\def\eg{e.g.}
\def\conj/{conjugate}
\def\deco/{decomposition}
\def\eg/{e.g.}
\def\fct/{function}
\def\gp/{group}
\def\hw/{highest weight}
\def\hwv/{highest weight vector}
\def\hwvs/{highest weight vectors}
\def\lw/{lowest weight}
\def\lwv/{lowest weight vector}
\def\lwvs/{lowest weight vectors}
\def\hds/{holomorphic discrete series}
\def\iff/{if and only if}
\def\inv/{invariant}
\def\irrde/{irreducible decomposition}
\def\meas/{measure}
\def\transf/{transform}
\def\rep/{representation}
\def\resp/{respectively}
\def\inters/{intertwines}
\def\interg/{intertwining}
\def\meta/{metaplectic representation}
\def\qu/{quaternion}
\def\rep/{representation}
\def\symdom/{ symmetric domain}
\def\st/{such that}
\def\shd/{subhead}
\def\transf/{transform}
\def\wrt/{with respect to}

\def\ra{\rightarrow}

\def\Norm#1#2#3{\Vert#1\Vert^{#3}_{{#2}}}

\def\tr{\operatorname{tr}}

\begin{document}
\def\abstractname{Abstract}
\def\chrefname{References}

\title[Eichler-Shimura isomorphism
]{Eichler-Shimura isomorphism for complex
hyperbolic lattices
 }
\author{Inkang Kim and  Genkai Zhang}
\address{School of Mathematics,
KIAS, Heogiro 85, Dongdaemun-gu
Seoul, 130-722, Republic of Korea}
\email{inkang@kias.re.kr}
\address{Mathematical Sciences, Chalmers University of Technology and
Mathematical Sciences, G\"oteborg University, SE-412 96 G\"oteborg, Sweden}
\email{genkai@chalmers.se}

\thanks{Research partially supported by
STINT-NRF grant (2011-0031291). Research by G. Zhang is supported partially
 by the Swedish
Science Council (VR). I. Kim gratefully acknowledges the partial support
of  grant  (NRF-2014R1A2A2A01005574) and a warm support of
Chalmers University of Technology during his stay.
}
\begin{abstract}
We consider the
cohomology group
$H^1(\Gamma,  \rho)$
of a discrete subgroup $\Gamma\subset
G=SU(n, 1)$ and the symmetric
tensor representation $\rho$
on $S^m(\mathbb C^{n+1})$.
 We
give an elementary proof of the
Eichler-Shimura isomorphism that harmonic
forms $H^1(\Gamma\backslash G/K, \rho)$
are $(0, 1)$-forms
for the automorphic holomorphic bundle
induced by the representation
$S^m(\mathbb C^{n})$ of $K$.
\end{abstract}

\maketitle

\baselineskip 1.35pc

\section{Introduction}

Let $B$ be the unit ball in
$\mathbb C^n
$ considered
as the Hermitian symmetric
space $B=G/K$ of $G=SU(n, 1)$, $n>1$. Let
$\Gamma$ be a cocompact torsion free discrete
subgroup of $G$ and $\rho$
a finite dimensional
representation of $G$, and $X=\Gamma\backslash B$.
The representation $\rho$ of $G$
defines also one for $\Gamma\subset G$.
The first cohomology  $H^1(\Gamma,  \rho)$
is of substancial interests
and appears naturally  in the study
of infinitesimal deformation of $\Gamma$
in a bigger group $G'\supset G$;
see \cite{Klingler-inv, KKP, GM}.
It is a classical result of
Raghunathan \cite{Rag-ajm} that the cohomology
group $H^1(\Gamma,  \rho)$
vanishes except when $\rho=\rho_m$
is the symmetric tensor $S^m(
\mathbb C^{n+1})$  (or $\rho_m^\prime$ on
$S^m(\mathbb C^{n+1})'$).
In a recent work \cite{Klingler-inv}
it is proved
that realizing   $H^1(X, \rho)$
as harmonic forms, it consists
of  $(0,1)$ forms
for the symmetric tensor of
the holomorphic tangent bundle
of $X=\Gamma\backslash B$.
The proof
in \cite{Klingler-inv}  uses a Hodge
vanishing theorem and the Koszul complex.
In the present paper we shall
give a rather elementary
proof of
the result.
We will prove that
any harmonic form with values
in $S^{m}(\mathbb C^{n+1})$
is $(0,1)$-form taking values
in $S^{m}(\mathbb C^n)$.
Let $TX$ and $T'X$ be the holomorphic tangent
and cotangent bundles respectively.
Let  $\mathcal L^{-1}$
be the
line bundle on $X$
defined
 so that $\mathcal L^{-(n+1)}$
is the canonical line bundle $\mathcal K=K_X$.
More precisely we shall prove the following, the notations being
explained in \S2,
\begin{theo+}
Let $\Gamma$ be a torsion free cocompact
lattice
of $G$ acting properly discontinuously on $B$.
\begin{enumerate}
\item
Let $\alpha\in A^1(\Gamma, B, \rho_m)$
be a harmonic form.
Then $\alpha$ is a $(0, 1)$ form
on $\Gamma\backslash B$ with values
in the 
holomorphic vector bundle 
 $S^m TX\otimes \mathcal L^{-m}$.
\item
Let $\alpha\in A^1(\Gamma, B, \rho_m^\prime)$
be a harmonic form.
Then $\alpha$ is a $(1, 0)$-form
on $\Gamma\backslash B$ with values
in the 
holomorphic vector bundle 
$S^m T'X \otimes \mathcal L^{-m}$
and $\alpha$ is symmetric in all $m+1$
variables. In particular $\alpha$
is naturally identified with a section of the
bundle $S^{m+1} T'X\otimes \mathcal L^{m}$.
\end{enumerate}
\end{theo+}

\begin{coro+}Let  $\Gamma$ be as  above
and assume that $\Gamma\backslash B$ is compact then
we have
$$
H^1(\Gamma, \rho_m)=
H^1(\Gamma\backslash B, S^m TX\otimes \mathcal L^{-m}),
\quad
H^1(\Gamma, \rho_m^\prime)=H^0
(\Gamma\backslash B, S^{m+1} T'X\otimes \mathcal L^m),
$$
where the cohomology on the right hand side
are the Dolbeault cohomology
of $\bar \partial$-closed $(0, 1)$
forms of the holomorphic vector bundles.
\end{coro+}

The case $n=1$, namely a Riemann surface
$\Gamma\backslash B$, is slightly
different. In that case the group cohomology
$H^1(\Gamma, \rho_{2j})$
 of the $2j$-th power
of the defining representation
of $\Gamma\subset SU(1, 1)$ will have
both holomorphic and antiholomorphic
components, $H^{(1, 0)}
(\Gamma, \rho_{2j})$,
 $H^{(0, 1)}(\Gamma, \rho_{2j})$, the
holomorphic part $H^{(1, 0)}(\Gamma, \rho_{2j})$
 corresponds
to
$$
H^{(1, 0)}(\Gamma, \rho_{2j})=H^{(1, 0)}(\Gamma\backslash B,
\mathcal K^{j+1}
)=H^{0}(\Gamma\backslash B,
\mathcal K^{j+1})$$
 of the tensor
power of the canonical line bundle. This
is known as the Eichler-Shimura correspondence;
see \cite[TH\'{E}OR\`{E}ME 1]{Shimura}
where a concrete construction was given.
We can also follow our proof and
get an  elementary proof of this result.

Our proof is  a bit tricky but
it is still very akin to the variation of
Hodge structures; conceptually
we are treating explicitly the filtration
of holomorphic bundles
defined by the central action of
$K$.
It is stated in \cite{Klingler-inv}
that the  results can be derived from the work of
  Deligne
and Zucker \cite{Zucker1, Zucker2}.
We note here that results of this type
that $(0, q)$-forms in the group cohomology
 $H^q(\Gamma, B, \rho)$
are actually $(0, q)$-forms for a corresponding
automorphic bundle have been obtained
much earlier by Matsushima
and Murakami \cite{MM1, MM2}. It seems that
one can prove the above result
by combining the works
of \cite{MM1, MM2, VZ, BW}.
But our method is down-to-earth hence we expect that it can be applied to various situations. For example we deal
with $n=1$ case, i.e., surface case in the last section, which is not available in \cite{Klingler-inv}.
 We will  investigate
further applications in a near future.

\section{Preliminaries}\label{pre}
Let $V=\mathbb C^{n+1}
$ be equipped with
the Hermitian inner product $\langle Jv, v\rangle$ of
signature $(n, 1)$, where
$J$ is the diagonal matrix $J=\diag( 1, \cdots, 1, -1)$
and $\langle v, w\rangle=\sum \bar v_i w_i$ the Euclidean form in $\mathbb C^{n+1}$.
We write $V=V_1\oplus\mathbb Ce_{n+1}$
with $V_1$ being the Euclidean space $\mathbb C^n$
with an orthonormal basis $\{e_k, k=1, \dots, n\}$.
Let $G=SU(n, 1)$ be the group
of linear transformations on $V$ preserving
the Hermitian form. The maximal compact subgroup
of $G$ is
$$
K=\{
\begin{bmatrix} A & 0\\
0& e^{i\theta}\end{bmatrix};
A \in U(n), \, e^{i\theta} \det A=1
\}=U(n),
$$
namely  $K=S(U(n)\times U(1))=U(n)$.
The subgroup $SU(n)\subset U(n)$
 viewed as a subgroup  in $K$ will be denoted by $SU(n)\times e$
to avoid confusion.
The Lie algebra
$\fg=\fsu(n, 1)$
 consists
of matrices $X$ such that
$X^\ast J + JX=0$.
The symmetric
space $G/K$ can be realized
as the unit ball $B$ in $V_1=\mathbb C^n$,
$B=G/K$ with $x_0=0$ being the base point.
Let $\fg=\fk\oplus \fp$ be the Cartan decomposition
of $\fg$ and the subspace $\fp=\{\xi_v; v\in \mathbb C^n\}$
with
$$
\xi_{v}=\begin{pmatrix} 0& v\\
\bar v& 0
\end{pmatrix}.
$$
The tangent space $T_{x_0}(B)$ at $x_0$
will be  identified with $\fp=\mathbb C^n$
as real spaces.

We fix an element in the center of the maximal compact
subalgebra $\fk=\fu(n)$
$$
H_0=(n+1)^{-1}\sqrt{-1} \text{diag}(1, \cdots, 1, -n),
$$
which defines the complex structure on $B$,
and we have
$$
\fsl(n+1)=\fsl(n) \oplus \mathbb C H_0
\oplus \fp^+ \oplus \fp^-.
$$
Then the holomorphic
and anti-holomorphic tangent
space $\fp^{\pm}$
consists of upper triangular, respectively
lower triangular matrices.
We denote
\begin{equation}
  \label{eq:xi-v}
\xi^+_v
=
\frac 12(
\xi_v
-i\xi_{iv})
=\begin{pmatrix} 0& v\\
0& 0
\end{pmatrix}
\in \fp^+,
  \,
\xi^-_v
=
\frac 12(
\xi_v
+i\xi_{iv})
=\begin{pmatrix} 0& 0\\
\bar v& 0
\end{pmatrix}\in \fp^-,
\end{equation}
the $\mathbb C$- and
 $\overline{\mathbb C}$-linear components of $\xi_v$.

Let $V_1=\mathbb C^n$ be the defining
representation
and
$\det(A) $ the
determinant representation
of $U(n)$.
We take   the diagonal
elements as  Cartan algebra
of $\mathfrak{gl}(n, \mathbb C)$
and the upper triangular matrices as positive
root vectors.
Denote
$\omega_1, \cdots, \omega_{n-1}$
the fundamental representations
of $U(n)$,
 so that $\omega_1=V_1$
is the defining representation above
and $\omega_{n-1}$ the dual representation. Note that $\omega_i$ has the highest weight $L_1+\cdots+L_i$ where
$L_j(\text{diag}(h_1,\cdots,h_n))=h_j$ is a canonical dual element on the Cartan algebra.

As complex representation of $\fu(n)$
we have
$$
\fp^+= \omega_1\otimes  \det
=V_1\otimes  \det, \,
\fp^-= \omega_{n-1}\otimes  \det^{-1}.
$$
This entails that, for $A\in U(n)$,
$$ A (\xi^+_{v_1}\wedge \cdots\wedge \xi^+_{v_n})=(\det A)^n A\xi^+_{v_1}\wedge \cdots\wedge A\xi^+_{v_n}=(\det A)^{n+1}
(\xi^+_{v_1}\wedge \cdots \wedge \xi^+_{v_n})
.$$
Hence
\begin{equation}\label{det}
K_X^{-1}=\wedge^n \fp^+=(\det)^{n+1}
\end{equation} and $\mathcal L=\det$.

We shall just identify
$\fp^+$ with $V_1$,
$\fp^+=V_1$, when the center action of
 $U(n)$ is irrelevant.

The defining representation $V$ of $G$
under $\fu(n)$ is
$$
V= V_1
 \oplus  \det^{-1}
$$
We shall consider its symmetric representation
$(S^{m}(V), \rho_m)$
of $G$
and $\fg$.
Note that we
have
\begin{equation}
  \label{eq:deco-2}
W=S^{m}(V)=\oplus_{k=0}^m W_k
=\oplus_{k=0}^m S^k(V_1)
\otimes e_{n+1}^{m-k},
\end{equation}
and we make the identification of the spaces
$$
W_k=S^k(V_1)
\otimes e_{n+1}^{m-k}
=S^k(V_1)
$$
whenever the factor  $e_{n+1}^{m-k} $ is irrelevant.

Note that the Euclidean inner product on $V$ induces
one on $W=S^{m}(V)$ and the above decomposition
is an orthogonal decomposition. Note also
that action of  $\rho_m(X)$
is Hermitian   for $X\in \fp$
and skew Hermitian for  $X\in \fk$.

A representation of $G$
on a finite dimensional real or
complex vector space
defines also a  vector bundle
 over the quotient space $\Gamma\backslash B$
and we recall briefly its construction
following
the exposition \cite{Rag-book, MM1} and also some notations there.
Let $(W, \rho)$ be a finite
dimensional representation of $G$ on
a real (or complex) vector space $W$. Eventually we shall only
consider $W=S^m(V)$
as above and its dual $S^m(V')$.
We fix
on $W$ a positive definite inner
(respectively) Hermitian  product 
so that $K$ acts  as orthogonally (resp. unitarily).
Let $\Gamma$ be a torsion free
discrete subgroup
of $G$. The restriction
of $\rho$ on $\Gamma$
will also be written as $\rho$.
 Suppose $\Gamma$ acts properly discontinuously
on $B$.
Let $\Gamma\times K
$ acts on $G\times W$
by
$(\gamma, \kappa)(g, w):=
(\gamma g \kappa^{-1}, \rho(\gamma)w)
$.
Then $E_\rho=G\times W/\Gamma\times K$
is a  vector bundle on $\Gamma\backslash B$.
The de Rham operator $d$ is well-defined
on $E_\rho$ and we let $\Delta_\rho=dd^\ast + d^\ast d$
be the corresponding Hodge Laplacian operator on
space of $p$-forms
 $\Omega(\Gamma\backslash B, E_\rho)$. We choose its standard realizations
as $W$-valued $p$-forms on $ G$ as follows.
Let $A^p(\Gamma, B, \rho)
$ be the space of $W$-valued
$p$-forms $\alpha$
on $G$ satisfying
\begin{description}
\item[(a)]
$\alpha(\gamma g)=\alpha(g)
$, $\gamma\in \Gamma$.
\item[(b)]
$
\rho(\kappa)\alpha(g\kappa^{-1})=
\alpha(g), \quad \kappa\in K.
$
\item[(c)]
$\iota(Y)\alpha=0$,  $Y\in \fk$.
\end{description}
Here $\iota(Y)$ is the pairing of
$Y\in \fg$ as left-invariant vector fields on $G$
(by differentiation
from right) with a $p$-form $\alpha$
on $G$, $\iota(Y)\alpha(Z_1, \cdots, Z_{p-1})
=\alpha(Y, Z_1, \cdots, Z_{p-1})$.
Equivalently it can be realized as $p$-forms
on $\Gamma\backslash G$ satisfying $(b)-(c)$ above and $A^p_0(\Gamma,B,\rho)$ denotes the space of $W$-valued $p$-forms on $\Gamma\backslash G$.
With some abuse of notation we denote $\Delta_\rho$
the corresponding Hodge Laplacian on
$A_0^p(\Gamma, B, \rho)$.

We shall also need the automorphic bundle defined by
representations of $K$, see \cite{MM1}. So let $(V, \tau)$ be a
complex
representation
of the complexification of $K_{\mathbb C}$
and we fix as above a Hermitian inner product
on $V$ so that $K$ acts unitarily.  The group
$\Gamma\times K$ acts on $G\times V$
by $(\gamma, \kappa)(g, w)=
(\gamma g\kappa^{-1}, \tau(\kappa)w)$.
Then $\mathcal E_{\rho}=\Gamma\times K \backslash
G\times V
$
defines a holomorphic vector bundle over $\Gamma\backslash B$.
The $p$-forms on the vector bundle
can be realized as the space $\mathcal A^p(\Gamma, B, \tau)$
(again with some abuse of notation)
of $p$-forms on $\Gamma\backslash G$ satisfying
\begin{description}
\item[(b')]
$
\tau(\kappa)\alpha(g\kappa^{-1})=
\alpha(g), \quad \kappa\in K.
$
\item[(c')]
$\iota(Y)\alpha=0$, $Y\in \fk$.
\end{description}

When $\rho$ is a complex representation of $G$ and $(\tau, K)$ is a sub-representation of $\rho$ restricted to $K$, then
 discrete group cohomology $H^{p+q}(\Gamma, B, \rho)$ and automorphic cohomology $H^{(p,q)}(\Gamma, B, \tau)$ are related by the work of \cite{MM1}.
\section{The Eichler-Shimura isomorphism}


In general, a real linear map $B$ on a complex vector space $W$ decomposes into $ \mathbb C$-linear part $B^+$ and $\overline{ \mathbb C}$-linear part $B^-$ so that $B(w)=B^+(w)+ B^-(\bar w)$ for $w\in W$.
For any real linear map $A: \fp\to \text{ End}_{\mathbb R}(W)$ from $\fp$
to  any complex vector
space $W$ we let
$$
A^+(\xi_v)
=
\frac 12(
A(\xi_v)
-iA(\xi_{iv})),\quad
A^-(\xi_v)
=
\frac 12(
A(\xi_v)
+iA(\xi_{iv}))
$$
be the $\mathbb C$-linear and respectively $\overline{\mathbb C}$-linear components.
In particular for any complex representation
$(W, \rho)$ of $G$ and $\fg$ we have
$$
\rho^{\pm}(\xi_v)=\rho(\xi^{\pm}_v),
$$
where $\xi_v^{\pm }$ are defined in (\ref{eq:xi-v}).
Let now $\rho=\rho_m$ be
the representation $S^m(V)$
and $\rho^m$ the dual representation
$S^m(V^\prime)$ of $\fg$. Note that $\rho_1$ is a defining representation $V$.
We start now with a few simple observations formulated
only $\rho=\rho_m$; the corresponding ones
hold for  $\rho^m$.

Denote by
$$P_k: W\to
W_k= S^k(V_1)\otimes e_{n+1}^{m-k}
$$ the orthogonal projection
onto the  component $W_k$
in  (\ref{eq:deco-2}), and write
$$
\alpha=\sum_{k=0}^m \alpha_k
$$
the corresponding decomposition for  $\alpha\in W=\sum_{k=0}^m W_k$.

Let $\{X_j\}$
 be an orthogonal basis
of $\fp$ viewed as tangent vectors
on $\Gamma\backslash G$ at a fixed
point $\Gamma g$ and $\{e_j\}$ be the corresponding orthonormal basis of $V_1$.
Let $T=T_\rho$
and $T^\ast=T_\rho^\ast$
be the operator defined on $A^1(\Gamma, B, \rho)$
as follows.
$$
T\alpha(X_1, X_{2})
= \rho(X_1)\alpha(X_2) -
\rho(X_2)\alpha(X_1)
$$
$$
T^\ast \alpha
=\sum_{j=1}^n \rho(X_j)\alpha(X_j)
$$

We recall the following result \cite[Corollary 7.50]{Rag-book}
\begin{prop+}Suppose $\alpha\in A_0^1(\Gamma, B, \rho)$
is harmonic, $\Delta_\rho \alpha=0$. Then
$T_\rho\alpha=0$ and $T_\rho^\ast\alpha=0$.
\end{prop+}

This can be restated as the following
(which is also proved  in \cite{KKP}
for $S^2(V)$ by using matrix computations).

\begin{coro+} Suppose $\alpha\in A_0^1(\Gamma, B, \rho)$
satisfies $T_\rho\alpha=0$ and $T_\rho^\ast\alpha=0$.
Then the $W$-valued $\mathbb R$-bilinear form
$(X, Y)\mapsto \rho(X)\alpha(Y)$ is symmetric
\begin{equation}
  \label{eq:symm-con}
\rho(\xi_v)\alpha(\xi_u)=
\rho(\xi_u)\alpha(\xi_v),
\end{equation}
and trace free
\begin{equation}
  \label{eq:tr-fr}
\sum_{j}
(
\rho(\xi_{e_j})\alpha(\xi_{e_j})
 +
\rho(\xi_{ie_j})\alpha(\xi_{ie_j}))
=0.
\end{equation}
\end{coro+}

Our theorem will be an easy consequence of the following proposition,
whose proof is based on a few elementary
lemmas.
\begin{prop+}\label{3.3}
\begin{enumerate}
\item
Suppose $\alpha\in {\text Hom}_{\mathbb R}(\fp, W)
$ satisfies
$T_\rho\alpha=T_\rho^\ast \alpha=0$.
Then $\alpha$ is $\overline{ \mathbb C}$-linear
and takes value in $W_m=S^mV_1$,
that is,
$\alpha=\alpha_m
=\alpha_m^-\in {\text Hom}_{\overline {\mathbb C}}(\fp^{-}, W_m)$.
\item  Suppose $\alpha\in {\text Hom}_{\mathbb R}(\fp, W^\prime)$
satisfies $T_{\rho^\prime}\alpha=T_{\rho^\prime}^\ast \alpha=0$.
Then  $\alpha$ is ${ \mathbb C}$-linear
and takes value in $S^m(V_1^\prime)$.
Moreover as an element in $
(\fp^{+})'\otimes S^m(V_1^\prime)
=(V_1)'\otimes S^m(V_1^\prime)
$,
it is symmetric in all variables, i.e.,
an element in
$ S^{m+1}(V_1^\prime)$,
 the leading
component
in $(V_1)'\otimes S^m(V_1^\prime)$.
\end{enumerate}
\end{prop+}

Denote $u^i v^{j-i}$
the symmetric tensor  power of $u$ and $v$ normalized
by
$$(u+v)^j=\otimes^j (u+v)=\sum_{i=0}^j \binom ji
u^i v^{j-i}.
$$

Note that the representation $\rho=\rho_m$ 
is the symmetric tensor $S^m(\bc^{n+1})$
and  $\rho'=\rho_m'$  its dual
throughout the paper.
\begin{lemm+}
\label{obs1}
\begin{enumerate}
\item
 Let $1\le k\le m-1$.
Then
for any $0\ne \xi_v\in \fp$,
$$
\rho(\xi_v): W_k\to W_{k+1}+W_{k-1},\,\,
\rho(\xi_v^+): W_k\to W_{k+1},\,\,
\rho(\xi_v^-): W_k\to W_{k-1},
$$
and on each  space it is nonzero.  Moreover
if $w\in W_k$ and $\rho(\xi_v^+) w=0$
or $\rho(\xi_v^-
) w=0$
for all $\xi_v^{\pm}\in \fp^{\pm}$ then $w=0$.
\item The restriction
 $\rho(\xi_v)|_{W_m}:
 W_m\to W_{m-1}$ on the top component
$W_m$ of $W$
is $\overline{\mathbb C}$-linear in $\xi_v$,
$\rho(\xi_v)|_{W_m}=\rho^-(\xi_v)|_{W_m}$,   and
$\rho(\xi_v)_{W_0}
$ on the bottom component
  is $ {\mathbb C}$-linear in $\xi_v$,
$\rho(\xi_v)_{W_0}=\rho^+(\xi_v)_{W_0}$.
\end{enumerate}
\end{lemm+}

\begin{proof}
The defining
representation $\rho_1$ is just
the matrix multiplication and we
have $$\rho_1(\xi_{v}) u=\langle v, u\rangle e_{n+1}$$
for $u\in V_1$, and
$$\rho_1(\xi_{v}) e_{n+1}= v.$$
Thus
$$\rho_1(\xi_{v}^+) u=0,\
\rho_1(\xi_{v}^-) u=\langle v, u\rangle e_{n+1},\
\rho_1(\xi_{v}) e_{n+1}= v,$$
and $$\rho_1(\xi_{v}^+) e_{n+1}=v,\
 \rho_1(\xi_{v}^-) e_{n+1}=0.$$ Taking the tensor power
we find
$$
\rho(\xi_{v}^+)
e_{n+1}^k=k v
e_{n+1}^{k-1}, \quad
\rho(\xi_{v}^-) e_{j}^k=k \overline{v_j}e_{n+1}
e_{j}^{k-1}, \quad 1\le j\le n,
$$
which are non-zero if $v_j\ne 0$.
Then
$$\rho(\xi^+_v)(W_k)=\rho(\xi^+_v)(S^k(V_1)\otimes e_{n+1}^{m-k})=(m-k)v e_{n+1}^{m-(k+1)}\in W_{k+1},$$$$
\rho(\xi^-_v)(W_k)=\rho(\xi^-_v)(S^k(V_1)\otimes e_{n+1}^{m-k})\in S^{k-1}(V_1)\otimes e_{n+1}^{m-(k-1)}\in W_{k-1}.
$$
First note that
$$
\rho(\kappa)\rho^{\pm}(\xi_v)\rho(\kappa^{-1})=
\rho^{\pm}(\xi_{\kappa v}), \quad \kappa\in SU(n)\times\{e\}, v\in V_1.
$$

If
$\rho(\xi_v^\pm) w=0$
for all $\xi_v^{\pm}\in \fp^{\pm}$
and for a fixed
$w\ne 0$,
then $$\rho(\kappa)\rho(\xi_v^\pm)\rho(\kappa^{-1})w=\rho(\xi_{\kappa v}^\pm)w=0$$ for all $\kappa\in SU(n)\times \{e\}$. Here the action of $K$ on $W$ is via the given representation $\rho$ from $G$. Hence it
is zero for all $\rho(\kappa^{-1})w$, and therefore zero
for $w=e_j^k, j=1,\cdots,n$, contradicting
 the previous claim.

The second part (2) follows immediately from
the above formulas for $\rho(\xi_v^{\pm})$
and fact that $W_{m+1}=0$ and 
 $W_{-1}=0$.
\end{proof}

The space $\text{Hom}_{\barc}(\fp^-, W_j)$
 of $\overline{\bc}$-linear forms $\beta=\beta^-$ on
$\fp^-$  will be identified
with  the tensor product
$(\fp^-)'  \otimes
W_j$. Using $(\fp^-)'=V_1\otimes \text{det}$,  the tensor product
is decomposed under  $K$ as \cite{Zlbk}
\begin{equation}
\begin{split}
\label{tensor-deco}
\text{Hom}_{\barc}(\fp^-, W_j)
&=(\fp^-)'  \otimes
S^j(V_1)  \otimes e_{n+1}^{m-j}\\
&\equiv (S^{j+1}(V_1) \otimes e_{n+1}^{m-j} ) \oplus
(S^{j-1, 1}(V_1)\otimes  e_{n+1}^{m-j})
\end{split}
\end{equation}
with the corresponding highest weights
$$
\omega_{1} \otimes j\omega_{1} =(j+1)\omega_{1}
 + ((j-1)\omega_{1}+\omega_{2}).
$$

\begin{lemm+}
\label{symm-vs-deco}
If
$\rho^-(\xi_u)\beta^-(\xi_v)=
\rho^-(\xi_v)\beta^-(\xi_u)$
then $\beta$ is in the first component
$S^{j+1}(V_1)$ in the above decomposition
(\ref{tensor-deco}).
\end{lemm+}
\begin{proof} Note that the relation
$\rho^-(\xi_u)\beta^-(\xi_v)=
\rho^-(\xi_v)\beta^-(\xi_u)$ is invariant under the $K$-action,
since
$$
\rho(\kappa)\rho^{\pm}(\xi_v)\rho(\kappa^{-1})=
\rho^{\pm}(\xi_{\kappa v}), \quad \kappa\in K, v\in V_1
$$ and
$$
\rho(\kappa)\beta(g\kappa^{-1})=\beta(g)
$$ for all $\kappa\in K$ (see Section \ref{pre}), which results in
$$\rho(\kappa)\rho^{\pm}(\xi_v)\beta(g\kappa^{-1})=\rho^{\pm}(\xi_{\kappa v})\beta(g).$$
 Thus
if $\beta^-$ satisfies the relation so is its component
in $((j-1)\omega_{1}+\omega_{2})$. We prove that any element
in $((j-1)\omega_{1}+\omega_{2})$
 satisfying the relation must
be zero. This space is an irreducible representation
of $K$ and we need only to check the relation for its highest weight vector.
The  highest weight vector of
$((j-1)\omega_{1}+\omega_{2})$ in $V_1\otimes S^j(V_1)
$
is
$$
\beta=\epsilon_2\otimes e_1^j - \epsilon_1\otimes (e_1^{j-1}e_2)
$$ where $\epsilon_i$ is a dual vector to $\xi^-_{e_i}$ in $\fp^-$.
We check the relation
$$
\rho(\xi^-_{e_2})\beta(\xi^-_{e_1})=
\rho(\xi^-_{e_1})\beta(\xi^-_{e_2}).
$$
The left hand side is $-e_1^{j-1}e_{n+1}$ whereas the
right hand side is $je_1^{j-1}e_{n+1}$, and the relation is
not satisfied. Hence $\beta$ should be in the first component $S^{j+1}(V_1)$.
\end{proof}
Note that $\{e_k\}$ is an orthogonal basis of $V_1$.
Observe that
for any $
\beta\in
\text{Hom}_{\barc}(\fp^-, W_j)$
we have
$$
\rho(\xi_v^+)\beta
\in \text{Hom}_{\barc}(\fp^-, W_{j+1}).
$$
\begin{lemm+}\label{3.7} Suppose $1\le j <m$.
The map
$$
T:
\text{Hom}_{\barc}(\fp^-, W_j)
\equiv (j+1)\omega_1
 \oplus
((j-1)\omega_1 +\omega_2)
\to W_{j+1}, \quad
\beta
\mapsto
\sum_{k=1}^n \rho(\xi_{e_k}^+) \beta(\xi_{e_k}^-)
$$
is up to non-zero constant an isometry on the space $
(j+1)\omega_1 $ where $\omega_i's$ are fundamental representations of $U(n)$ introduced in Section \ref{pre}.
\end{lemm+}
\begin{proof}
It is clear that $T$
is a $K$-intertwining map
from $\text{Hom}_{\barc}(\fp^-, S^j(V_1))$
 into $W_{j+1}$. By Schur's
lemma it is either zero or an isometry
up to non-zero constant
on the irreducible space
$(j+1)\omega_1 $.
To find the constant we take
$\beta=
\varepsilon_1
\otimes e_1^j e_{n+1}^{m-j}
$
where $\varepsilon_1$ is the dual form of $\xi_{e_1}^-$. It
is indeed in the first component
$ (j+1)\omega_1 $ and is actually the highest weight
vector.
Then
by direct computation we find
$$
T\beta=(m-j)e_1^{j+1}
e_{n+1}^{m-j-1},
$$
which is nonzero.
\end{proof}

We consider the corresponding symmetry property
for the dual representation $\rho'=\rho^m$.
\begin{lemm+}
\label{sym-pro}
Suppose  $\beta=\beta^+$ is  $S^m(V_1^\prime)$-valued
${\mathbb C}$-linear
 form on $\fp^+$.
If
 $
\rho'(\xi_u^+)
\beta(\xi_v^+)=\rho'(\xi_v^+)
\beta(\xi_u^+)$
 then $\beta$ as an element
in $(\fp^+)'  \otimes
S^m(V_1^\prime) $ is symmetric in all $
m+1$ variables.
\end{lemm+}
\begin{proof}
The statement is equivalent to
that
$\beta(\xi_v^+)(
\xi_{u}^+, \xi_{v_1}^+,
\cdots, \xi_{v_{m-1}}^+)$
is symmetric in  all
$m+1$ variables.
However the equality
$\rho'(\xi_u^+)
\beta(\xi_v^+)=\rho'(\xi_v^+)
\beta(\xi_u^+)$  implies that it
is symmetric in the first two variables and thus
is symmetric in all $m+1$ variables. 
More precisely, viewing $\rho'(\xi^+_u)\beta(\xi^+_v)$ 
and $\rho'(\xi^+_v)\beta(\xi^+_u)$ as  elements in $S^m(V')$,
$$\rho'(\xi^+_u)\beta(\xi^+_v)(e_{n+1},\cdots,e_{n+1})=\beta(\xi^+_v)(\rho'(\xi^+_u)
e_{n+1},\cdots, 
\rho'(\xi^+_u) e_{n+1})$$
$$= \rho'(\xi^+_v)\beta(\xi^+_u)(e_{n+1},\cdots,e_{n+1})=
\beta(\xi^+_u)(\rho'(\xi^+_v) e_{n+1},\cdots, \rho'(\xi^+_v) e_{n+1}).$$ Hence from
$\rho'(\xi^+_u) e_{n+1}=u$ and $\rho'(\xi^+_v) e_{n+1}=v$ 
and identifying $\fp^+=V_1$, we get
$$\beta(\xi^+_v)(\xi^+_u,\cdots,\xi^+_u)=\beta(\xi^+_u)(\xi^+_v,\cdots,\xi^+_v).$$
\end{proof}

We prove now Proposition 3.3.
\begin{proof}
 We shall prove  by  induction
 that all $\alpha_j=0$ for $k\le m-1$.
Let $1\le k \le m-1$. Taking the $k$-th component of (\ref{eq:symm-con})
we get
\begin{equation}
\label{kth-eq-1}
 \rho^{+}(\xi_v) \alpha^{+}_{k- 1}(\xi_u)
=
\rho^{+}(\xi_u)
\alpha^{+}_{k-1}(\xi_v),
\end{equation}
\begin{equation}
\label{kth-eq-2}
\rho^{-}(\xi_u) \alpha^{-}_{k+ 1}(\xi_v)
=
\rho^{-}(\xi_v)
\alpha^{-}_{k+ 1}(\xi_u),
\end{equation}
\begin{equation}
\label{kth-eq-3}
\rho^+(\xi_u)
\alpha^-_{k-1}(\xi_v)
=
\rho^-(\xi_v)
\alpha^+_{k+1}(\xi_u).
\end{equation}

We prove first that $\alpha_0=0$.
Consider the $1$-component of the identity
\begin{equation}
  \label{eq:trace-0}
T_\rho^\ast \alpha=
\sum_j \left (
\rho(\xi_{e_j})\alpha(\xi_{e_j})
+
\rho(\xi_{ie_j})\alpha(\xi_{ie_j})
\right)=0
\end{equation}
and write each term in terms of their $\mathbb C$-linear
and  $\overline{\mathbb C}$-linear parts. Note
that bilinear $\mathbb C$-linear and
bilinear $\overline{\mathbb C}$-linear terms have their sum zero.
Also  on the component $W_0$ the action $\rho(\xi_u)
=\rho(\xi_u^+)$ is $\mathbb C$-linear, by Lemma \ref{obs1}.
Thus
$$
\sum_j
\left (
\rho^+(\xi_{e_j})
\alpha_0^-(\xi_{e_j})
+
\rho^-(\xi_{e_j})
\alpha_2^+(\xi_{e_j})
\right)=0.
$$
But by the equality of  (\ref{kth-eq-3})
for $k=1$ we have $\rho(\xi_{e_j}^-)
\alpha_2^+(\xi_{e_j})=\rho(\xi_{e_j}^+)
\alpha_0^-(\xi_{e_j})$.
Namely
\begin{equation}
  \label{eq:pr-1}
2\sum_j
\rho(\xi_{e_j}^+)
\alpha_0^-(\xi_{e_j})
=0.
\end{equation}
Taking inner product with $e_1 e_{n+1}^{m-1}\in W_1$,
and using the fact that
$$
\langle
\rho(\xi_{e_1}^+)
\alpha_0^-(\xi_{e_1}), e_1 e_{n+1}^{m-1}\rangle
=\langle
\alpha_0^-(\xi_{e_1}),
\rho(\xi_{e_1}^-)(e_1 e_{n+1}^{m-1})
\rangle
=\langle
\alpha_0^-(\xi_{e_1}),
e_{n+1}^{m}
\rangle
$$
and
$$
\langle
\rho(\xi_{e_j}^+)
\alpha_0^-(\xi_{e_j}), e_1 e_{n+1}^{m-1}\rangle
=\langle
\alpha_0^-(\xi_{e_j}),
\rho(\xi_{e_j}^-)(e_1 e_{n+1}^{m-1})
\rangle
=0, j\ne 1,
$$
 we see
that $\langle \alpha_0^-(\xi_{e_1}), e_{n+1}^{m}\rangle=0
$, namely $\alpha_0^-(\xi_{e_1})=0
$ since it is a scalar multiple of 
$ e_{n+1}^{m}$. By the $K$-invariance of above relation (3.8)
we may replace $e_1
$ by any $e_j$, and get  $\alpha_0^-(\xi_{e_j})=0$,
i.e.,  $\alpha_0^-=0$ and $\alpha_0$ is $\mathbb C$-linear,
$\alpha_0=\alpha_0^+$. Now $W_0=\mathbb Ce_{n+1}^m$
is one-dimensional and $\alpha_0$
is thus of the form $$\alpha_0(\xi_u)=\alpha_0(\xi^+_u)=\langle u_0, u\rangle
e_{n+1}^m
$$ for some $u_0\in V_1$.
The relation
  (\ref{kth-eq-1} )
 implies that
$$
\langle u_0, u\rangle v e_{n+1}^{m-1}
=\langle u_0, v\rangle u e_{n+1}^{m-1}
$$
for all $u, v\in V_1$. This is impossible
unless $u_0=0$ since $\dim V_1>1$, i.e.,
$\alpha_0=0$.

Taking  the $0$-th component of the equality $
\rho(\xi_u)
\alpha(\xi_v)=\rho(\xi_v)
\alpha(\xi_u)
$ we get
$$
\rho^-(\xi_u)
\alpha_1(\xi_v)=\rho^-(\xi_v)
\alpha_1(\xi_u).
$$
Changing $v$ to $iv$ we find
$$
\rho^-(\xi_u)
\alpha_1(\xi_{iv})=-i\rho^-(\xi_v)
\alpha_1(\xi_u).
$$
Summing the two results we get
$$
\rho^-(\xi_u)(\alpha_1(\xi_{iv})
+i\alpha_1(\xi_{v})
)
=0.
$$
Taking further the inner product with $e_{n+1}^m\in W_0$
we have
\begin{equation*}
\begin{split}
0&=\langle (\rho^-(\xi_u)(\alpha_1(\xi_{iv})
+i\alpha_1(\xi_{v}), e_{n+1}^m\rangle
=\langle \alpha_1(\xi_{iv})
+i\alpha_1(\xi_{v}), \rho^+(\xi_u)
e_{n+1}^m\rangle
\\
&
=\langle ((\alpha_1(\xi_{iv})
+i\alpha_1(\xi_{v}),
u e_{n+1}^{m-1}\rangle
\end{split}
\end{equation*}
for all $u$. Thus
$\alpha_1(\xi_{iv})
+i\alpha_1(\xi_{v})=0$, namely $\alpha_1$ is $\bar \bc$-linear,
$\alpha_1=\alpha_1^-$. Furthermore it follows from Lemma
\ref{symm-vs-deco}
that $\alpha_1$ is an element in the component
$S^2(V_1)$ in $(\mathfrak p^-)' \otimes S^1(V_1)$.

We take now   the $2$-component of the identity (\ref{eq:trace-0})
using again the fact that $\alpha_1$ is $\barc$-linear,
and find
$$
0=\sum_j \left (
\rho^+(\xi_{e_j})\alpha_1(\xi_{e_j})
+
\rho^+(\xi_{ie_j})\alpha_1(\xi_{ie_j})
\right)
=2
\sum_j \left (
\rho^+(\xi_{e_j})\alpha_1(\xi_{e_j}^-)
\right).
$$
But $\alpha_1$ is in the component
 $2\omega_1=
S^2(V_1)$ and Lemma
\ref{3.7}
implies that $\alpha_1=0$.

Using the above procedure successively we prove
then that $\alpha_{j}=0$ for $j\le m-2$.
Consequently we have $\alpha_{m-1}^+ =0$
and $\alpha_{m-1} =\alpha_{m-1}^{-}$.
Taking the trace of $(m-2)$-th component
of   (\ref{eq:tr-fr})
we have again
$\sum_j \rho^+(\xi_{e_j})\alpha_{m-1}^-(\xi_{e_j})=0$
and $\alpha_{m-1}=0$ by the same arguments.

Finally we consider  the $(m-1)$-th component of the equality $
\rho(\xi_u)
\alpha(\xi_v)=\rho(\xi_v)
\alpha(\xi_u)
$. We have
$$
\rho^-(\xi_u)
\alpha_m(\xi_v)=\rho^-(\xi_v)
\alpha_m(\xi_u).
$$
Replacing $u$ by $iu$  gives
$$
-i\rho^-(\xi_u)
\alpha_m(\xi_v)=\rho^-(\xi_v)
\alpha_m(\xi_{iu}).
$$
Thus
$$\rho^-(\xi_v)\alpha_m^+(\xi_u)
=\frac 12
\rho^-(\xi_v)\left(
\alpha_m(\xi_u)-i\alpha_m(\xi_{iu})\right)
=0.
$$
This holds for all $\xi_v\in \fp$.
Thus  $\alpha_m^+(\xi_u)=0$ by Lemma \ref{obs1},
and $\alpha_m$
is $\barc$-linear.

The second part (2) of the Proposition on the dual representation
can be proved similarly using the similar arguments
and
Lemma \ref{sym-pro}.
\end{proof}

We prove now  Theorem 1.1 and Corollary 1.2.

\begin{proof} The statements in Theorem 1.1
follows from Proposition 3.3. Indeed if $\alpha\in
A^1(\Gamma, B, \rho_m)$
is a harmonic form,
then it will have values in $S^m(\bc^n)$ by Proposition \ref{3.3}.
By the relation $\bc^n=\fp^+\otimes \det^{-1}$ we have
$$S^m(\bc^n)=(\fp^+)^m\otimes (\det)^{-m}=S^m TX\otimes \mathcal L^{-m},$$
proving that $\alpha$ is a $(0, 1)$-section of
$S^m TX\otimes \mathcal L^{-m}$.
 The proof of the second one is similar.
The claim that $\alpha$ is $\barc$-linear
is precisely that $\alpha$ is a $(0, 1)$-form. This
proves the first part, and the second part follows similarly
from Proposition 3.3 (2).

Let  $\alpha$
be a harmonic form representing an element  $ H^1(\Gamma, \rho)$.
Write $\alpha=\sum_{k=0}^m \alpha_k$
according to the decomposition (\ref{eq:deco-2}). It
follows then from
above  that $\alpha_k=0$ for $k<m$,
i.e. $\alpha=\alpha_m$.
The isomorphism
of the cohomology $H^1(\Gamma, \rho)$ and
$H^1(\Gamma\backslash B, S^m TX\otimes \mathcal L^{-m})$
is  then a consequence of
\cite[Proposition 4.2 and Theorem 6.1]{MM1}.
The second isomorphism is proved similarly.
\end{proof}

\section{The Eichler-Shimura isomorphism
for Riemann surfaces and applications}

We consider now the case $n=1$. Keeping
the previous notation we consider the group
cohomology $H^1(\Gamma, \rho_m)$ of
the tensor power $S^m(\mathbb C^2)$
of the representation  of $\Gamma\subset SU(1, 1)$.
In this case $H^1(\Gamma, \rho_m)$
has
a decomposition as
$H^1(\Gamma, \rho_m)=H^{(1, 0)}(\Gamma, \rho_m)
+H^{(0, 1)}(\Gamma, \rho_m)
$, and in contrast
to the case $n\ge 2$ the
component $H^{(1, 0)}(\Gamma, \rho_m) $
is not vanishing but it is dual to
$H^{(0, 1)}(\Gamma, \rho_m)$.
This Eichler-Shimura isomorphism further
gives  a correspondence
between $H^{(1, 0)}(\Gamma, \rho_m) $
and $H^0$-cohomology of a line bundle over the
Riemann surface $\Sigma:=\Gamma\backslash B$.
We denote $ K_\Sigma=\mathcal K$
the holomorphic
cotangent
bundle, i.e. the canonical line bundle
on $\Gamma\backslash B$.
\begin{theo+}  Realizing
 $H^{1}(\Gamma, \rho_m) $
as the space of harmonic forms
on $\Gamma\backslash B$
we have  $
H^1(\Gamma, \rho_m)=H^{(1, 0)}(\Gamma, \rho_m)
+H^{(0, 1)}(\Gamma, \rho_m)
$ and furthermore the two space are dual to each other,
$$
H^{(1, 0)}(\Gamma, \rho_m)
=H^0(\Sigma, \mathcal K^{\frac m2 +1}), \quad
H^{(0, 1)}(\Gamma, \rho_m)=H^0(\Sigma, \mathcal K^{\frac m2 +1})^\ast.
$$
\end{theo+}
\begin{proof}
 We prove the second isomorphism
using the computation in \S3.
Let $\alpha$ be a $(0, 1)$
form in $H^{(0, 1)}(\Gamma, \rho_m)$. Using
 $z\in B$ near $z=0$ as local coordinate
as above, let
$\alpha=\sum_{j=0}^m \alpha_j $ be the decomposition
of $\alpha  $ in the decomposition of $S^m(\bc^2)
=\oplus_{j=0}^m\mathbb C e_1^j e_2^{m-j}$. The symmetry condition
(3.6) implies that
$$
\rho^+(\xi_u)\alpha(\xi_v)=
\rho^+(\xi_u)\alpha^-(\xi_v)=
\rho^-(\xi_v) \alpha^+(\xi_u)=0
$$
since $\alpha$ is $\overline{\mathbb C}$-linear,
hence $\alpha=\alpha^-$ and $\alpha^+(\xi_u)=0$. Thus
$$(\rho^+(\xi_u)\alpha(\xi_v))_j=
\rho^+(\xi_u)
\alpha(\xi_v)_{j-1}
=0$$
for all $j\ge 1$. But then since
$\rho^+(\xi_{e_1})$ maps $e_1^{j-1}e_2^{m-j+1}$
to $e_1^{j}e_2^{m-j}$ for $1\le j\le m $,
the component $\alpha(\xi_v)_{j-1}$
is vanishing for all $1\le j\le m$,
and we have $\alpha=\alpha_m$. Now  from equation (\ref{det}), $K_\Sigma^{-1}=\det^2$. Hence  $S^m V_1=S^mT\Sigma\otimes \det^{-m}=\mathcal K^{-m}\otimes \mathcal K^{\frac{m}{2}}=\mathcal K^{-\frac m2 }$
and we have thus
$\alpha\in H^1(\Sigma, \mathcal K^{-\frac m2 })$
which is dual to
$H^0(\Sigma, \mathcal K^{\frac m2  +1})$ by Serre duality. That
this map is onto is a consequence of the general
results of \cite{MM1, MM2}, as in the proof of Corollary 1.2.
\end{proof}

We give now an application of the above result
computing the tangent space of
the Hitchin's Teichm\"uller component
of representations
of $\Gamma$ in a semisimple Lie
group $G$.
We shall only treat the case
$G=SL(n, \mathbb R)$ even though
much computations can be carried
over to other cases.  The result might be known to experts
but it seems still to provide some novel understanding
for the geometry of the component.

We consider two representations of $SL(2, \mathbb R)=SU(1, 1)$
into the group $SL(n, \mathbb R)$
and compute the corresponding group
cohomologies  of $\Gamma\subset SL(2, \mathbb R)$.
We consider first the real
representation of
 $SL(2, \mathbb R)$
 on the symmetric tensor
$(S^m(\mathbb R^2), \rho_m)
$  in the group
 $SL(m+1, \mathbb R)$.
Let $\tau_k$, $k\le m$, be the
representation $\rho_k$ in $
SL(k+1, \mathbb R)$ considered
as a representation in $SL(m+1, \mathbb R)$.
We compute
the corresponding cohomologies
which can be viewed as the tangent
space of the variety at the
respective points.

\begin{theo+}\label{Zariski} Realizing
  the elements in the group cohomologies as harmonic forms we have that
$H^1(\Gamma, \rho_m, \mathfrak{sl}(m+1, \mathbb R))$ and
$H^1(\Gamma, \tau_k, \mathfrak{sl}(m+1, \mathbb R))$ are
real forms in the
space
$$
\sum_{j=1}^{m}
H^0(\Sigma,
\mathcal K^{j+1}) +H^0(\Sigma,
\mathcal K^{j+1})^\ast
$$
and
$$
\sum_{j=1}^k
H^0(\Sigma,
\mathcal K^{j+1}) +H^0(\Sigma,
\mathcal K^{j+1})^\ast +
(H^0(\Sigma,
\mathcal K^{\frac k2+1})
+H^0(\Sigma,
\mathcal K^{\frac k2+1})^\ast)^{2(m-k)}$$$$
+(H^0(\Sigma, \mathcal K) +H^0(\Sigma,\mathcal K)^\ast)^{(m-k)^2}.
$$
\end{theo+}
\begin{proof} We consider
the complexification
$\mathfrak{sl}(2, \mathbb C)$
and its representation
$S^k\mathbb C^2$
in $\mathfrak{sl}(m+1, \mathbb C)$.
The real representation
$S^k\mathbb R^2$
of $\mathfrak{sl}(2, \mathbb R)$
in
$\mathfrak{sl}(m+1, \mathbb R)$
is the fixed point of the conjugation
$X\to \bar X$ of
$\mathfrak{sl}(m+1, \mathbb C)$. Now the
adjoint representation of
$\mathfrak{sl}(2, \mathbb C)$ under $\rho_m$ in
$\mathfrak{sl}(m+1, \mathbb C)$
 is decomposed as
\cite{Kostant, Hitchin}
$$
\mathfrak{sl}(m+1, \mathbb C)
=\sum_{j=1}^{m} S^{2j}\mathbb C^2
$$
with the first component
$S^{2}\mathbb C^2$
being $\mathfrak{sl}(2, \mathbb C)$ itself.
Now by Theorem 4.1 we have
$$
H^{(1, 0)}(\Gamma, S^{2j}\mathbb C^2)
=H^{(1,0)}(\Sigma, \mathcal K^{j})=
H^{0}(\Sigma, \mathcal K^{j+1}),$$$$
\quad
H^{(0, 1)}(\Gamma, S^{2j}\mathbb C^2)
=H^{(0, 1)}(\Sigma, \mathcal K^{-j})
=H^{0}(\Sigma, \mathcal K^{j+1})^\ast.
$$
The involution $X\mapsto \bar X$ on the
one-forms is now $f(z)(dz)^{j} \otimes  dz
\mapsto  {\bar f(z)} \partial_z^{j} \otimes  d\bar z$.
Thus the real cohomology
$H^1(\Gamma, \mathfrak{sl}(m+1, \mathbb R))$
is a real form in the space stated.
Now  under  the action $\tau_k$
we have
$$
\mathfrak{sl}(m+1, \mathbb C)
=\sum_{j=1}^k S^{2j}\mathbb C^2
\oplus (S^{k}\mathbb C^2)^{2(m-k)}
\oplus \mathbb C^{(m-k)^2}
$$
the cohomology of $\Gamma$
in the $S^{2j}\mathbb C^2$
is computed as above. The cohomology
in $\mathbb C
$ is
$$
H^1(\Gamma, \mathbb C)
=H^0(\Sigma, \mathcal K)
+H^0(\Sigma, \mathcal K)^\ast
$$
the space of abelian differentials. The rest
of the claim follows immediately.
\end{proof}

The set $$\{\rho_m\circ \phi|\phi:\Gamma\rightarrow SL(2,\mathbb R)\ \text{discrete and faithful}\}/\sim$$ constitutes the Fuchsian locus
in the Hitchin component. The above theorem shows that the tangent space of Hitchin component at Fuchsian locus consists of $\sum_{j=2}^{m+1}H^0(\Sigma, K_\Sigma^j)$. When $m=2$, the Hitchin component is the set of convex real projective structures on a surface. Furthermore it is known that the Hitchin component is a holomorphic vector bundle over Teichm\"uller space with fibers cubic holomorphic forms \cite{La,Lo}. In the forthcoming paper \cite{KZ}, we will analyze this case in more details to show that the Hitchin component is a K\"ahler manifold.

\end{document}